\documentclass[draft, a4paper]{article}
\usepackage{amsmath, amsfonts, amsthm, amssymb}
\numberwithin{equation}{section}

\newcommand{\Q}{\mathbb{Q}}

\newcommand{\N}{\mathbb{N}}
\newcommand{\Z}{\mathbb{Z}}
\renewcommand{\leq}{\leqslant}
\renewcommand{\geq}{\geqslant}

\newcommand{\cQ}{\overline{\Q}}
\newcommand{\Proj}{\mathbb{P}}

\newcommand{\ep}{\varepsilon}
\newtheorem{lemma}{Lemma}

\newtheorem*{thm}{Theorem}

\theoremstyle{definition}
\newtheorem*{ack}{Acknowledgements}

\renewcommand{\u}{\mathbf{u}}
\newcommand{\lab}{\label}
\newcommand{\bfP}{\mathbb{P}}

\newcommand{\h}{\mathrm{h.c.f.}}
\newcommand{\ma}{\mathbf}
\newcommand{\beq}{\begin{equation}}
\newcommand{\eeq}{\end{equation}}
\newcommand{\ve}{\varepsilon}
\newcommand{\mcal}{\mathcal}
\newcommand{\la}{\lambda}

\newcommand{\x}{{\bf x}}

\newcommand{\bxi}{\boldsymbol{\xi}}

\begin{document}

\title{Simultaneous equal sums of three powers}

\author{T.D. Browning$^1$ and D.R. Heath-Brown$^2$\\
\small{$^1$\emph{School of Mathematics,  
Bristol University, Bristol BS8 1TW}}\\
\small{$^2$\emph{Mathematical Institute,
24--29 St. Giles',
Oxford OX1 3LB}}\\
\small{$^1$t.d.browning@bristol.ac.uk}, \small{$^2$rhb@maths.ox.ac.uk}}
\date{}
\maketitle

\begin{abstract}
Using a result of Salberger \cite{s} we show that the number of
non-trivial positive integer solutions $x_0,\ldots,x_5 \leq B$ to the
simultaneous equations
$$
x_{0}^{c}+x_{1}^{c}+x_{2}^{c}=x_{3}^{c}+x_{4}^{c}+x_{5}^{c},
\quad
x_{0}^{d}+x_{1}^{d}+x_{2}^{d}=x_{3}^{d}+x_{4}^{d}+x_{5}^{d},
$$
is $o(B^3)$ whenever $d>\max\{2,c\}$.
\end{abstract}

\section{Introduction}

The purpose of this short note is to apply recent work of Salberger to
an old problem in analytic number theory.
More precisely we shall see how Salberger's new results \cite{s} 
concerning the distribution of rational points on projective algebraic
varieties can be used to study the number of
positive integer solutions to the simultaneous equations
\begin{equation}\lab{2six}
x_{0}^{c}+x_{1}^{c}+x_{2}^{c}=x_{3}^{c}+x_{4}^{c}+x_{5}^{c},
\quad
x_{0}^{d}+x_{1}^{d}+x_{2}^{d}=x_{3}^{d}+x_{4}^{d}+x_{5}^{d},
\end{equation}
in a region $\max x_{i}\leq B$, for fixed positive integers $c<d$.
There are clearly $6B^{3}+O(B^{2})$
trivial solutions in which $x_{3},x_{4},x_{5}$
are a permutation of $x_{0},x_{1},x_{2}$.  We write $\mathcal{N}_{c,d}(B)$ for
the number of non-trivial solutions, and our primary goal is to
estimate this quantity.  

When $c=1$ and $d=2$ it is rather easy
to show that 
$$
\mathcal{N}_{1,2}(B)= CB^3\log B(1+o(1)),
$$ 
for an
appropriate constant $C>0$, so that the non-trivial solutions
dominate the trivial ones.   In all other
cases we would like to know that the trivial solutions dominate
the non-trivial ones.  This has only been established when 
$c=1$, or when
$c=2$ and $d=3$ or $4$. Specifically, 
it has been shown by Greaves \cite{greaves} that
\beq\lab{prev1}
\mathcal{N}_{1,d}(B)\ll_{d,\ep} B^{\frac{17}{6}+\ep},
\eeq
and by Skinner and Wooley \cite{s-w} that
\beq\lab{prev2}
\mathcal{N}_{1,d}(B)\ll_{d,\ep} B^{\frac{8}{3}+\frac{1}{d-1}+\ep}.
\eeq
This latter result reproduces Greaves' result for $d=7$, and improves upon it for
$d\geq 8$.  Moreover, work of  Wooley \cite{wooley} shows that
\beq\lab{prev3}
\mathcal{N}_{2,3}(B)\ll_{\ep} B^{\frac{7}{3}+\ep}, 
\eeq
and Tsui and Wooley \cite{t-w} have shown that
$$
\mathcal{N}_{2,4}(B)\ll_{\ep} B^{\frac{36}{13}+\ep}.
$$
We are now ready to record the contribution that we have been able to
make to this subject.

\begin{thm}
Let $\ep>0$, and suppose that $c,d$ are positive integers such that
$c<d$ and $d\geq 4$. 
Then we have
$$
\mathcal{N}_{c,d}(B)\ll_{c,d, \ep} 
B^{\frac{11}{4}+\ep} +
B^{\frac{5}{2}+\frac{5}{3cd}+\ep}.
$$
\end{thm}

The implied constant in this estimate is allowed to depend at most
upon $c,d$ and the choice of $\ve$.  When $c=1$, 
it is easy to see that our result improves upon 
\eqref{prev1} for $d \geq 6$, and upon 
\eqref{prev2} for $d \leq 11$.  Moreover, it retrieves
\eqref{prev1} for $d=5$.  When $c<d$ are arbitrary positive integers
such that $d \geq 4$, it follows from the theorem that
$$
\mathcal{N}_{c,d}(B)=o(B^3).
$$
An application of Greaves' bound \eqref{prev1}
shows that the same is true when $c=1$ and $d=3$.  Similarly, the
bound \eqref{prev3} of Wooley handles the case $c=2$ and $d=3$.
This therefore confirms the paucity of non-trivial solutions to the
pair of equations \eqref{2six} for all positive integers $c,d$ such
that $d>\max\{2,c\}$. 

A crucial aspect of our work involves working in
projective space.  It is clearly natural to talk about rational points
on varieties, rather than integral 
solutions to systems of equations.  Whereas a single integer solution to
\eqref{2six} can be used to generate infinitely many others by scalar
multiplication, all of these will actually correspond to the same
projective rational point on the variety defined by \eqref{2six}.  
The transition between estimating $\mcal{N}_{c,d}(B)$ and counting 
non-trivial projective rational points of bounded height will be made
precise in the following section.

For distinct positive
integers $c,d$ the pair of equations \eqref{2six} defines a projective algebraic
variety $X_{c,d}\subset \bfP^5$ of dimension $3$. The trivial
solutions then correspond to rational points lying on certain planes contained in
$X_{c,d}$.  Our proof of the theorem makes crucial use of a rather
general result due to Salberger \cite{s}, that provides a good upper bound
for the number of rational points of bounded height that lie on the
Zariski open subset formed by deleting all of the planes from an
arbitrary threefold in $\bfP^5$.  
In order to apply this result effectively we shall need to study
the intrinsic geometry of $X_{c,d}$. This will be carried out
separately in the final section of this paper.

Our primary goal in this paper was merely to obtain exponents strictly
less than 3.  However, in private communications with the authors, Salberger has indicated
how the upper bound in our theorem can be substantially improved upon
by applying work of his that is currently in preparation. The  
authors are grateful to him for this observation, and a
number of useful comments that he made about an earlier version of
this paper.

\begin{ack}
Part of this work was undertaken while the first author was attending
the \emph{Diophantine Geometry} intensive research period at the
Centro di Ricerca Matematica Ennio De Giorgi in Pisa, the 
hospitality of which is gratefully acknowledged.  
While working on this paper, the first
author was also supported at Oxford University by
EPSRC grant number GR/R93155/01.  
\end{ack}

\section{Proof of the theorem}\lab{final}

Let $X_{c,d}\subset \bfP^5$ denote the
threefold defined by the pair of equations \eqref{2six}, for 
positive integers $c<d$ such that $d \geq 4$. Our first task is to consider the possible 
planes contained in $X_{c,d}$, for which it will clearly suffice to
consider the planes contained in $F_d$, where $F_k \subset \bfP^5$
denotes the non-singular Fermat hypersurface
\beq\lab{fermat}
x_{0}^{k}+x_{1}^{k}+x_2^k+x_3^k+x_4^k+x_{5}^{k}=0,
\eeq
for any $k\in \N$.
Now there are certain obvious planes contained in $F_{k}$, which we refer to as
``standard''. These planes are obtained by first partitioning the
indices $\{0,\ldots,5\}$ into three distinct pairs.  
For each such pair  $\{i,j\}$, one then associates a vector
$\x_{i,j}=(x_i,x_j) \in (\overline{\Q}\setminus\{0\})^2$ such that
$$
x_i^k+x_j^k=0.
$$
Then if $J_0,J_1,J_2$ are the sets of indices so formed, and 
$[\x_{J_0},\x_{J_1},\x_{J_2}] \in \bfP^{5}$ denotes the
corresponding point in $F_k$, we thereby obtain the plane
$$
\{[\la_0\x_{J_0},\la_1 \x_{J_1}, \la_2 \x_{J_2}]:
~[\la_0,\la_1,\la_2]\in \bfP^2\} \subset F_k.
$$
It is not hard to see that  this procedure produces exactly $15 k^{3}$ standard planes
contained in $F_{k}$. The following result shows that these are the only planes
contained in $F_{k}$, if $k$ is at least $4$.

\begin{lemma}\lab{standard}
Let $k\geq 4$. Then any plane contained in $F_{k}$ is standard.
\end{lemma}

It ought to be remarked that for $k \geq 5$ the classification of
planes in $F_k$ is a straightforward consequence of the well-known classification of
lines in $F_k$, as discussed by Debarre \cite[\S 2.5]{debarre}, for example.
Thus the chief novelty of Lemma~\ref{standard} is that we are also able to
handle the case $k=4$.  Lemma \ref{standard} will be established  in the next section.
An important consequence of this result is that
any plane in the threefold $X_{c,d}$ must correspond to a standard
plane in $F_d$, since we have assumed that $d\geq 4$.
Among the planes contained in $X_{c,d}$ are the six ``trivial'' planes
$$
x_0=x_i, \quad 
x_1=x_j, \quad x_2=x_k,
$$
where $\{i,j,k\}$ is a permutation of the set $\{3,4,5\}$.

We are now ready to complete the proof of the theorem.
In estimating
$\mathcal{N}_{c,d}(B)$, it is clearly enough to count
primitive vectors $\x=(x_0,\ldots,x_5)\in \N^6$, where $\x$ is said to
be ``primitive'' if $\h(x_0,\ldots,x_5)=1$.  
In general, if we know that there are $O(B^{\delta})$
primitive vectors that make a contribution to $\mathcal{N}_{c,d}(B)$,
for some constant $\delta>1$, then there will be  
$$
\ll \sum_{m\leq B}(B/m)^{\delta}\ll B^{\delta}
$$ 
vectors in total.  
For any rational point $x=[\x] \in \bfP^5(\Q)$ such that
$\x \in \Z^{6}$ is primitive, we shall write 
$$
H(x)=\max_{0\leq i \leq 5}|x_i|
$$
for its height.  
For any $B \geq 1$ and any Zariski open subset $U\subseteq
X_{c,d}$ we shall set 
$$
N_U(B)=\#\{x \in U\cap\bfP^5(\Q): ~H(x) \leq B\}.
$$
A little thought reveals that the trivial planes are the only planes
contained in $X_{c,d}$ that can possibly contain rational points
$[x_0,\ldots,x_5]\in \bfP^5(\Q)$ in which $x_0,\ldots,x_5$ share the same sign.
Since we are only interested in counting non-trivial
solutions to the pair of equations \eqref{2six}, 
it will therefore suffice to estimate $N_{U}(B)$ in the case that 
$U\subset X_{c,d}$ is the Zariski open subset formed by deleting all of the planes from 
$X_{c,d}$.

It is now time to record the result that forms the backbone of our
estimate for $N_{U}(B)$. 
The following upper bound is due to Salberger \cite[Remark 8.6]{s}.

\begin{lemma}\lab{salb}
Let $\ep>0$, and suppose that $Z\subset\Proj^{5}$ is a geometrically integral 
threefold defined over $\overline{\Q}$, of degree $D$.  Let 
$Y$ be the complement of the union of all planes contained in $Z$.  Then we have
$$
\#\{x\in Y\cap \bfP^5(\Q): ~H(x)\leq B\} \ll_{D,\ep} B^{\frac{11}{4}+\ep} +
B^{\frac{5}{2}+\frac{5}{3D}+\ep}.
$$
\end{lemma}

The proof of Lemma \ref{salb} follows from applying a birational projection
argument to a corresponding bound for threefolds in $\bfP^4$.  This
latter result is provided by combining work of the authors'
\cite[Theorem 3]{n-2} with the proof of Salberger's earlier estimate 
\cite[Theorem 3.4]{s@ENS} for hypersurfaces in $\bfP^4$.

Before applying Lemma \ref{salb} to our situation it is clear that we
shall need to say something about the integrality and degree of the
threefold $X_{c,d}$, as given by \eqref{2six}.  
In fact we shall show in Lemma \ref{d-i}, in the next section, 
that $X_{c,d}$ is geometrically integral and has degree $cd$.
Subject to the verification of these facts below, we may therefore
conclude from Lemma \ref{salb} that 
$$
N_{U}(B)\ll_{c,d,\ep}
B^{\frac{11}{4}+\ep} +
B^{\frac{5}{2}+\frac{5}{3cd}+\ep},
$$
for any $\ep>0$.  This completes the proof of the theorem.

\section{Geometry of Fermat varieties}

Throughout this section let $c,d$ be positive integers such that 
$c<d$, and let $X_{c,d}\subset \bfP^5$ denote the threefold \eqref{2six}.
Our goal is to study the geometry of $X_{c,d}$, with a view to
establishing the facts that were employed in the previous section.  

Let us begin by considering the singular locus of $X_{c,d}$. Now if
$\bxi=[\xi_0,\ldots,\xi_5]$ is a singular point of $X_{c,d}$ with at
least two non-zero coordinates $\xi_i,\xi_j$, then it follows from the
Jacobian criterion that $\xi_i^{c-1}\xi_j^{d-1}=\delta
\xi_i^{d-1}\xi_j^{c-1}$, for some $\delta \in \{-1,+1\}$.  
Hence the $2(d-c)$-th powers of any two non-zero
coordinates of $\bxi$ must coincide, and so we may conclude that the singular locus of
$X_{c,d}$ is finite.  We are now ready to establish the following
result.

\begin{lemma}\lab{d-i}
$X_{c,d}$ is a geometrically integral variety of degree
$cd$.
\end{lemma}

The first part of Lemma \ref{d-i} may be compared with work of Kontogeorgis
\cite[Theorem 1.2]{kont}, who has established that the variety
$\sum_{i=0}^n x_i^c=\sum_{i=0}^n x_i^d=0$ is geometrically reduced and
irreducible if $n\geq 5$.

In order to prove Lemma \ref{d-i},  we define
$$
G_k: \quad x_{0}^{k}+x_{1}^{k}+x_2^k-x_3^k-x_4^k-x_{5}^{k}=0,
$$
in analogy to \eqref{fermat}.  Then $G_k$ is a non-singular
hypersurface of dimension $4$, and it is clear that $X_{c,d}=G_{c}\cap
G_{d}$ is a complete intersection in $\bfP^5$.  Hence it follows
from \cite[Proposition II.8.23]{hart} that $X_{c,d}$ is
Cohen-Macaulay as a subscheme of $G_d$.  
Since the singular locus of
$X_{c,d}$ is finite, and so has codimension $\geq 2$ in $X_{c,d}$, we may
therefore apply \cite[Theorem 18.15]{eisen} to deduce that $X_{c,d}$ is
geometrically reduced and irreducible.  Turning to the degree of
$X_{c,d}$, it is not hard to check that $G_c$ and $G_d$ intersect
transversely at a generic point of $X_{c,d}$, again using the fact
that $X_{c,d}$ has finite singular locus.  But then an
application of B\'ezout's theorem, in the form
\cite[Theorem 18.3]{harris}, immediately reveals that 
$$
\deg X_{c,d}=\deg G_{c}.\deg G_{d}=cd.
$$
This completes the proof of Lemma \ref{d-i}.


Our final task in this section is to establish Lemma \ref{standard}.
Let $k \geq 4,$ and suppose that  we are given a plane $\Pi$ contained
in $F_k$, as given by \eqref{fermat}.
Then $\Pi$ must be generated by three non-collinear points $e_0,e_1,e_2 \in
F_k$.  We may assume after a linear change of variables that 
$$
e_0=[1,0,0,\ma{e}_0], \quad 
e_1=[0,1,0,\ma{e}_{1}], \quad 
e_2=[0,0,1,\ma{e}_{2}], 
$$
for certain $\ma{e}_i=(e_{i,3},e_{i,4},e_{i,5})\in \cQ^3$.
But then there exist linear forms of the shape
$$
L_0=u_0, \quad L_1=u_1, \quad L_2=u_2, \quad L_i=a_iu_0+b_iu_1+c_iu_2,
$$
for $i=3,4,5$, that are defined over $\cQ$ and satisfy
\beq\lab{viva''}
L_0(\u)^k+L_1(\u)^k+L_2(\u)^k+L_3(\u)^k+L_4(\u)^k+L_5(\u)^k=0,
\eeq
identically in $\u=(u_0,u_1,u_2)$.  In particular we may henceforth assume that
none of the forms $L_3,L_4,L_5$ are identically zero, since it is
well-known that a non-singular hypersurface of dimension $3$ contains
no planes.
We proceed to differentiate the identity \eqref{viva''} with respect to $u_0$, giving
\begin{equation}\lab{viva;}
u_0^{k-1} + a_3L_3(\u)^{k-1}+a_4L_4(\u)^{k-1}+a_5L_5(\u)^{k-1}=0.
\end{equation}
On writing $M_{i}(\u)= a_i^{\frac{1}{k-1}}L_i(\u)$ for
$i=3,4,5$, we must therefore investigate the possibility that 
\beq\lab{viva}
u_0^{k-1}+M_3(\u)^{k-1}+M_4(\u)^{k-1}+M_5(\u)^{k-1}=0,
\eeq
identically in $\u$.  Now either 
there exist $\la_3,\la_4,\la_5\in \cQ$, not all
zero, such that $M_i(\u)=\la_i u_0$ for $i=3,4,5$, 
or we have produced a projective linear space of positive dimension that is
contained in the surface $y_0^{k-1}+y_1^{k-1}+y_2^{k-1}+y_3^{k-1}=0$.  
Since $k\geq 4$, the only such
spaces contained in this surface are the obvious lines.  Thus we may
assume, without loss of generality, 
that in either case there exist constants $\la_3,\la_4 \in \cQ$ such that 
$$
M_3(\u)=\la_3 u_0, \quad M_4(\u)=\la_4 M_5(\u),
$$
with $1+\la_3^{k-1}=0$ and $(1+\la_4^{k-1})M_4(\u)=0$ identically.
We claim that the only possibility here is that $a_4=a_5=0$ in the
definition of $M_4,M_5$.  To
see this it clearly suffices to suppose for a contradiction that $a_4a_5 \neq 0$, since
$L_4, L_5$ are non-zero in \eqref{viva;}.
But then, on considering the identity \eqref{viva''} satisfied by the 
original forms $L_0,\ldots,L_5 \in \cQ[\u]$, we may deduce that there are
constants $\alpha,\beta \in \cQ$ such that
$$
\alpha u_0^k+u_1^k+u_2^k=\beta L_5(\u)^k,
$$
identically in $\u$.  This is clearly impossible for $k \geq 4$, and
so establishes the claim.

In terms of the original linear forms $L_0,\ldots,L_5$
satisfying \eqref{viva''}, our consideration of \eqref{viva} has therefore led to the conclusion
$$
L_3(\u)=\mu_3 u_0, \quad a_4=a_5=0,
$$
for $\mu_3$, which must clearly be a $k$-th root of $-1$.
It remains to consider the possibility that there exist
$b_4,b_5,c_4,c_5 \in \cQ$ such that 
$$
u_1^k+u_2^k+(b_4u_1+c_4u_2)^k+(b_5u_1+c_5u_2)^k=0,
$$
identically in $u_1,u_2$.  But this corresponds to a line on the Fermat
surface of degree $k$, and so leads to the
conclusion that $L_4=\mu_4 L_i, L_5=\mu_5L_j$ for some permutation
$\{i,j\}$ of $\{1,2\}$, where $\mu_4,\mu_5$ are appropriate $k$-th
roots of $-1$.   This completes
the proof that any plane $\Pi$ contained in $F_k$ must be standard if
$k \geq 4$, as claimed in Lemma \ref{standard}.

Let $m \geq 1$ and $k \geq 3$ be integers.
It is interesting to remark that the argument just presented can be
easily generalised to treat possible $m$-dimensional 
linear spaces that are contained in the non-singular hypersurface
$$
y_0^k+\cdots+y_{2m+1}^k=0,
$$
in $\bfP^{2m+1}$.  Thus it can be shown that all such linear spaces
are the obvious ones, and that there are precisely
$c_m k^{m+1}$ of them, where
$$
c_m=(2m+1)\cdot (2m-1) \cdots 3 \cdot 1.
$$


\begin{thebibliography}{99}

\bibitem{n-2} T.D. Browning and D.R. Heath-Brown, Counting rational
points on hypersurfaces. {\em J. reine angew. Math.}, {\bf 584} (2005), 83--115.

\bibitem{debarre} O. Debarre,
{\em Higher-dimensional algebraic geometry}. Springer-Verlag, 2001.

\bibitem{eisen} D. Eisenbud, {\em Commutative algebra. With a view
  toward algebraic geometry}. Springer-Verlag, 1995. 


\bibitem{greaves}
G.R.H. Greaves, Some Diophantine equations with almost all 
solutions trivial. {\em Mathematika,} {\bf 44} (1997), 14--36. 

\bibitem{harris}
J. Harris, {\em Algebraic geometry.} Springer-Verlag,
1995.

\bibitem{hart}
R. Hartshorne, 
{\em Algebraic geometry}. Springer-Verlag, 1977.

\bibitem{kont}
A. Kontogeorgis, 
Automorphisms of Fermat-like varieties.
{\em Manuscripta Math.} {\bf 107} (2002), no. 2, 187--205.

\bibitem{s@ENS} P. Salberger, 
Counting rational points on hypersurfaces of low dimension.  {\em
  Ann. Sci. \'Ecole Norm. Sup.} {\bf 38}  (2005),  no. 1, 93--115.

\bibitem{s} P. Salberger, Rational points of bounded height on
  projective surfaces. {\em Submitted}, 2005.

\bibitem{s-w}
C.M. Skinner and T.D. Wooley, On the paucity of 
non-diagonal solutions in certain diagonal Diophantine systems.
{\em Quart. J. Math. Oxford Ser. (2),} {\bf 48} (1997), 255--277. 

\bibitem{t-w}
W.Y. Tsui and T.D. Wooley, The paucity problem for 
simultaneous quadratic and biquadratic equations. {\em
Math. Proc. Cambridge Philos. Soc.,} {\bf 126} (1999), 209--221.

\bibitem{wooley}
T.D. Wooley, An affine slicing approach to certain paucity problems,
{\em Analytic number theory, (Allerton Park, IL, 1995).} 803--815, 
Prog. Math., {\bf 139}, Birkh\"{a}user Boston, 1996. 

\end{thebibliography}
\end{document}